\documentclass[12pt]{amsart}
\usepackage{latexsym, amssymb, amsmath, amscd}
 \usepackage{eucal}
 \usepackage{graphics}

    \newtheorem{rema}{Remark}[section]
    \newtheorem{propo}[rema]{Proposition}

   \newtheorem{theo}[rema]{Theorem}
   \newtheorem{def-theo}[rema]{Definition-Theorem}
 
   \newtheorem{defi}[rema]{Definition}
    \newtheorem{lemma}[rema]{Lemma}
    \newtheorem{corol}[rema]{Corollary}

\newcommand{\Z}{\mathbb Z}

\newcommand{\Q}{\mathbb Q}
\newcommand{\N}{\mathbb N}
\newcommand{\C}{\mathbb C}
\newcommand{\R}{\mathbb R}
\newcommand{\A}{\mathbb A}
\newcommand{\p}{\partial}
\newcommand{\f}{\phi}

\newcommand{\al}{\alpha}

\newcommand{\EE}{\mathcal{E}}

\newcommand{\cA}{\mathcal{A}}
\newcommand{\im}{\operatorname{Im}}
\newcommand{\epfv}{\hspace{1em}$\Box$\vspace{1em}}

\begin{document}

\title[Locally Finite Derivations and $\mathcal{E}$-Derivations]{On  Images Of  Locally Finite Derivations and $\mathcal{E}$-Derivations}

\author{Arno van den Essen and Wenhua Zhao}

\date{August 09, 2022}

 \address{A. van den Essen, Department of Mathematics, Radboud University
Nijmegen, The Netherlands. Email: A.vandenEssen@math.ru.nl, arno.vd.essen@gmail.com}
\address{W. Zhao, Department of Mathematics, Illinois State University,
Normal, IL 61761. Email: wzhao@ilstu.edu}

\begin{abstract}  
Some cases of the LFED Conjecture, proposed by the second author \cite{Open-LFNED}, for certain integral domains 
are proved. In particular, the LFED Conjecture is completely established for the field of fractions $k(x)$ of the polynomial algebra $k[x]$, the formal power series algebra $k[[x]]$ and 
the Laurent formal power series algebra $k[[x]][x^{-1}]$, where $x=(x_1, x_2, \dots, x_n)$ denotes $n$ commutative free variables and $k$ a field of characteristic zero. Furthermore, the relation between the LFED Conjecture and the Duistermaat-van der Kallen Theorem \cite{DK} is also discussed and emphasized. 
\end{abstract}
  
   

\thanks{The second author has been partially supported 
by the Simons Foundation grant 278638}

 \bibliographystyle{alpha}
\maketitle

\section*{\bf Introduction}

The aim of this paper to draw some attention to and prove some cases for a recent conjecture of the second author, concerning images of locally finite derivations and $\mathcal{E}$-derivations. An $\mathcal{E}$-derivation is a map of the form $\delta\!:=1-\f$, where $\f$ is a ring endomorphism. Furthermore, if $\f$ is a $k$-algebra endomorphism of a $k$-algebra $\cA$, we say 
$\delta$ is a {\it $k\EE$-derivation} of $\cA$. 
There is an abundance of literature related to kernels of derivations. See for example \cite{N}, \cite{E} and \cite{F} and the references therein. The study of these kernels
has turned out to be very fruitful in the solution of several problems
in affine algebraic geometry.  See for example \cite{ML}, \cite{B}, \cite{D}, \cite{G}, \cite{K} and \cite{E}. 
However, images of derivations and $\EE$-derivations, by contrast, have hardly been investigated.

In the recent paper \cite{Open-LFNED} the second author proposed two conjectures related to images of locally finite or locally nilpotent derivations and $\EE$-derivations,
in the context of the so-called Mathieu subspaces or Mathieu-Zhao spaces (MZ-spaces for short). These spaces (see Definition \ref{Def1.3}) were introduced in \cite{Z1} and \cite{MS} to study
the Jacobian Conjecture and several related problems. 
In this paper we will consider one of these two conjectures for commutative algebras, 
the so-called LFED Conjecture. It asserts the following: {\it if $\cA$ is an algebra over a field $k$ of characteristic zero and $\delta$ a locally finite $k$-derivation or a locally finite $k\mathcal{E}$-derivation of $\cA$, then $\im\delta$ is an MZ-space of $\cA$}. For some other studies on 
this conjecture (and also on the other conjecture), 
see \cite{EWZ} and \cite{Open-LFNED}-\cite{OneVariableCase}.

A strong motivation for investigating this conjecture is given by a remarkable theorem obtained by Duistermaat and van der Kallen in \cite{DK}. To describe it let $\C[x,x^{-1}]\!:=\C[x_1,\dots,x_n, x_1^{-1},\dots,x_n^{-1}]$ be the  Laurent polynomial algebra in $n$ variables $x$ over $\C$. Each $f\in\C[x,x^{-1}]$ can be written in a unique way as a finite
sum $\sum_{a\in\Z^n} f_a x^a$, with $f_a\in\C$. The Duistermaat-van der Kallen theorem states that {\it the subspace 
of $f\in\C[x,x^{-1}]$ such that $f_0=0$ is an MZ-space of $\C[x,x^{-1}]$}. 

Since the proof of this theorem is highly non-trivial
and since this theorem is tightly connected with several conjectures that all imply the Jacobian Conjecture, it is very desirable
to have a new proof for it. It is at this point that the LFED Conjecture becomes very interesting. 
Namely, if the LFED Conjecture is true, then the truth of this conjecture implies the Duistermaat-van der Kallen theorem. 

To see how let $p_1,\dots, p_n$ be the first $n$ (distinct) 
prime numbers and define a $\C$-algebra endomorphism $\f$ of $\C[x,x^{-1}]$ by setting for all  
$1\leq i\leq n$,  
$$
\f(x_i)\!:=p_i x_i.
$$
\noindent One readily verifies that $\f$ is locally finite. Furthermore, set 
$\delta\!:=1-\f$ and $x^a\!:=x_1^{a_1}\dots x_n^{a_n}$,
where $a\!:=(a_1,\dots,a_n)\in\Z^n$.  Then $\f(x^a)=p_1^{a_1}\dots p_n^{a_n}x^a$. Consequently, we have 
$$
\im \delta=\{f\in\C[x,x^{-1}]\,\big|\, f_0=0\}.
$$
\noindent Therefore the Duistermaat-van der Kallen theorem follows if the LFED Conjecture holds for this $\C$-algebra endomorphism $\f$.  Similarly, if
$c_1,\dots, c_n\in\C$ are linearly independent over $\Q$ and $D\!:=\sum_i c_ix_i\p/\p x_i$, one easily verifies that $D$
is a locally finite $\C$-derivation on $\C[x,x^{-1}]$ such that $\im D$ is precisely the subspace  of $f\in\C[x,x^{-1}]$ with $f_0=0$.
So the derivation case of the LFED Conjecture also implies the Duistermaat-van der Kallen theorem.\\

{\bf Arrangement}:  In the first section we recall some definitions and  well-known
results concerning ($\EE$-)derivations and MZ-spaces. In section two we prove some cases of the LFED Conjecture for local integral domains over a field $k$ of characteristic zero. In particular, the LFED Conjecture is completely established in this section for the field of fractions $k(x)$ of the polynomial algebra $k[x]$ and the formal power series algebra $k[[x]]$. 
Using some of these results we show in the last section that the LFED Conjecture holds for the $k$-algebras $k[[x]][x^{-1}]$. So in particular we obtain local analogues of the Duistermaat-van der Kallen theorem.   
 
\renewcommand{\theequation}{\thesection.\arabic{equation}}
\renewcommand{\therema}{\thesection.\arabic{rema}}
\setcounter{equation}{0}
\setcounter{rema}{0}
\setcounter{section}{0}

\section{\bf Preliminaries}

Throughout this section $k$ is a field and $\cA$ a commutative $k$-algebra.
A {\em $k$-derivation} on $\cA$ is a $k$-linear map $D$ such that $D(ab)=D(a)b+aD(b)$ 
for all $a,b\in \cA$. Using induction on $n$ we get
that for all $n\geq 1$ and all $a,b\in \cA$,  
\begin{align}\label{Eq1.1}
D^n(ab)=\sum_{i=0}^n {n\choose i} D^i(a) D^{n-i}(b).
\end{align}

A $k$-derivation $D$ on $\cA$ is said to be {\em locally finite} if for each $a\in \cA$ the $k$-vector subspace spanned by the elements
$D^m(a)$ ($m\geq 0$) is finite dimensional over $k$.
A $k$-derivation $D$ is said to be {\em locally nilpotent} if for each $a\in \cA$ there exists an $m\geq 1$ such that
$D^m(a)=0$. For such a derivation we can define the {\em $D$-degree} of a non-zero element $a\in \cA$, denoted by $\deg_D a$, to be the smallest positive integer 
$n$ such that $D^{n+1}(a)=0$. The following result (\cite[Proposition 1.3.32 ii)]{E}) will play a crucial role in the next sections. To keep this paper self-contained we here include a short proof of it, which is given in \cite{F}.

\begin{propo}\label{Prop1.1}
Let $\cA$ be an integral domain containing $\Q$ and $D$ a locally nilpotent derivation
on $\cA$. If $D(a)=ba$, with $a,b\in \cA$ and $a\neq 0$, then $b=0$. 
\end{propo}

\medskip

\noindent{\bf Proof:}  Assume $b\neq 0$. Let $n=\deg_D a$ and $m=\deg_D b$. So $D^{n+1}(a)=0$ and $D^{m+1}(b)=0$.
It follows from Eq.\,(\ref{Eq1.1}) that $D^{n+m+1}(ab)=0$ and that $D^{n+m}(ab)={n+m\choose n}D^n(a)D^m(b)$. Since $D^n(a)$ and $D^m(b)$ are non-zero, $\cA$ is a domain and $\cA$ is a $\Q$-algebra, we get that
 $D^{n+m}(ab)\neq 0$. So $\deg_D ab=n+m$. Since $\deg_D D(a)=n-1$, the equality $D(a)=ab$ gives $n-1=n+m$, a contradiction.
\epfv

Besides the $k$-derivations described above we will also need the notion of a $k\EE$-derivation:

\begin{defi}\label{Def1.2}
 A $k\EE$-derivation on $\cA$ is a map of the form $\delta\!:=1_\cA-\f$, where $\f$ is
a $k$-algebra endomorphism of $\cA$. We will write $1$ instead of $1_\cA$ if the $k$-algebra $\cA$ is clear in the context. 
\end{defi}
 
As above a $k\EE$-derivation $\delta$ of $\cA$ is said to be {\em locally finite} if for each $a\in \cA$ the $k$-vector subspace spanned by the  elements $\delta^m(a)$, with $m\geq 1$, is finite dimensional over $k$    
and $\delta$ is said to be {\em locally nilpotent} if for each $a\in \cA$ there exists an $m\geq 1$ such that $\delta^m(a)=0$.

\medskip

In order to formulate the main conjecture, we recall the notion of a Mathieu-Zhao space. 
First, if $M$ is a $k$-subspace of $\cA$ we define its {\em radical}, 
denoted by $\mathfrak r(M)$, 
as the set of $a\in \cA$ such that $a^m\in M$ 
for all large $m$, i.e., for all $m$ greater than some $N$.

\begin{defi}\label{Def1.3}
 Let $M$ be a $k$-linear subspace of $\cA$. Then $M$ is called a Mathieu-Zhao space (MZ-space)
of $\cA$ if the following holds: if $a\in \mathfrak r(M)$ then for each $b\in \cA$ there exists $N_b\geq 1$ such that $ba^m\in M$ for all $m\geq N_b$. 
\end{defi}

Some remarks on the definition above are as follows. First, the notion was first introduced by the second author \cite{Z1} under the name {\em Mathieu subspaces}. 
It was later suggested by the first author \cite{E2} to change it 
to {\em Mathieu-Zhao spaces}. 
Second, the definition above is slightly different from 
the original definition \cite[Definition $4.1$]{Z1}. But it can be shown 
(e.g., see \cite[Proposition $2.1$]{MS}) that these two definitions are actually equivalent to each other.

Third, an ideal in $\cA$ is an example of an MZ-space of $\cA$, so an MZ-space is a generalization of the notion of an ideal in a ring.
In contrast to ideals it is often very hard to decide if a given subspace is an MZ-space or not. For example, it is still unknown if the kernel $\ker L$ is an MZ-space of $\cA$ for most of $k$-linear functionals $L: \cA\rightarrow k$ of $\cA$. On the other hand, as came to light by the work of the 
second author, various problems in affine algebraic geometry boil down to the question {\it whether or not for some specific $k$-algebras $\cA$ and specific $k$-linear functionals of $\cA$ the kernel $\ker L$ of $L$ is an MZ-space of $\cA$}. In particular, the Duistermaat-van der Kallen theorem mentioned in the
introduction is the statement that $\ker L$ is an MZ-space of $\C[x,x^{-1}]$, where
$L:\C[x,x^{-1}]\rightarrow\C$ is the $\C$-linear functional defined by $L(1)=1$ and $L(x^a)=0$ for all $a\in\Z^n\backslash\{0\}$.

\medskip

Now let us turn to the main conjecture (which was formulated in \cite{Open-LFNED}). 
From now on {\it $k$ will be a field
 of characteristic zero}.

\medskip

\noindent{\bf The LFED Conjecture.} {\em Let $D$ be a locally finite $k$-derivation (resp. $\delta$ a locally finite $k\EE$-derivation)
on $\cA$. Then $\im D$ (resp. $\im \delta$) is an MZ-space of $\cA$.}

\medskip

In the sequel we frequently use the following fact: if $\overline{k}$ is an algebraic closure of $k$ and the LFED Conjecture holds for
$\bar{\cA}:=\bar{k}\otimes_k\cA$, then it holds for $\cA$.

To see this extend the $k$-basis $\{1\}$ of $k$ to a $k$-basis $\{1\}\cup\{e_i \,|\, i\in I\}$ of $\bar{k}$. So every element
of $\bar{k}\otimes_k\cA$ can be written uniquely as $1\otimes c+\sum_i e_i\otimes c_i$ for some $c,c_i\in\cA$. Let $D$
be a locally finite $k$-derivation of $\cA$. Extend it to a $\bar{k}$-derivation $\bar{D}$ on $\bar{\cA}$. Then $\bar{D}$
is also locally finite. Assume $\im \bar{D}$ is an MZ-space of $\bar{\cA}$. To see that $\im D$ is an MZ-space of $\cA$, let
$a^m\in \im D$ for all large $m$. Then $(1\otimes a)^m=1\otimes a^m\in \im \bar{D}$ for all large $m$. So if $b\in\cA$, then
$1\otimes ba^m=(1\otimes b)(1\otimes a)^m\in \im \bar{D}$ for all large $m$. Hence for some $c,c_i\in \cA$ we have
$$1\otimes ba^m=\bar{D}(1\otimes c+\sum_i e_i\otimes c_i)=1\otimes Dc+\sum_i e_i\otimes Dc_i$$
It follows that $ba^m=Dc\in \im D$ for large $m$. So $\im D$ is an MZ-space of $\cA$. A similar argument works
for  $k\EE$-derivations $\delta$ of $\cA$.

As pointed out in \cite{Open-LFNED}, if $k$ is a field of characteristic $p>0$, both statements of the LFED Conjecture fail. To see this observe first that if
an MZ-space $M$ contains $1$, then $1^m\in M$ for all $m\geq 1$, which implies that $M=\cA$. Now take for example $D\!:=d/d t$ on the univariate polynomial algebra $k[t]$.
Then clearly $1\in \im D$, but $t^{p-1}\notin \im D$. So $\im D$ is not an MZ-space of $k[t]$. 

Furthermore, let
$\f\!:=\exp(D)$, i.e., $\f(t)=t+1$ and $\delta=1-\f$. Then $1=(1-\f)(-t)\in \im \delta$. 
However $t^{p-1}\notin \im  \delta$: for if $t^{p-1}=u(t)-u(t+1)$ for some $u(t)\in k[t]$, 
then replacing $t$ by $t+i$ in the equation for each $0\le i\le p-1$, and adding the resulted $p$ equations together we get 
$$
\sum_{i=0}^{p-1} (t+i)^{p-1}=0.  
$$
Since by Fermat's Little Theorem $i^{p-1}=1$ for all $1\le i\le p-1$, 
setting $t=0$ in the equation above we get $p-1=0$, a contradiction.  
  
\medskip

To conclude this section we recall Proposition $5.4$ of \cite{Open-LFNED} 
for commutative $k$-algebras and give a somewhat different proof.   
 
\begin{propo}\label{Prop1.4}
 Let $\cA$ be a commutative ring contained in a $\Q$-algebra and $\f$ a ring endomorphism of $\cA$ such that $\f^i=\f^j$ for some $1\leq i<j$. Put $M\!:=(1-\f)(\cA)$. Then $\mathfrak r(M)=\mathfrak r(I)$, where $I$ is the ideal of $a\in \cA$ such that $\f^r(a)=0$ for some $r\geq 1$. Furthermore, $M$ is an MZ-space of $\cA$.
\end{propo}

\noindent{\bf Proof:} Let $a\in I$, say $\f^r(a)=0$. Then, using   $a=((1-\f)+\f)^r(a)=0$ we get that $a\in (1-\f)(\cA)=M$.
So $I\subseteq M$. Hence $\mathfrak r(I)\subseteq \mathfrak r(M)$. To prove the converse observe that $0=\f^i-\f^j=\f^ig(\f)(1-\f)$, where $g(\f)=1+\f+\cdots+\f^{j-i-1}$. Hence $(1-\f)\cA\subseteq \ker \f^i g(\f)$. So if $a\in \mathfrak r(M)$, i.e., $a^m\in (1-\f)\cA$ for
all large $m$, then 
$$
\f^i(a)^m+(\f^{i+1}(a))^m+\cdots+(\f^{ j-1}(a))^m=0, \mbox{ for all large } m.
$$

\noindent It then follows from lemma \ref{lma1.5} below  that  $\f^i(a)$ is nilpotent. So $\f^i(a)^r=0$
for some $r\geq 1$, i.e., $\f^i(a^r)=0$. So $a\in \mathfrak r(I)$. Hence $\mathfrak r(M)=\mathfrak r(I)$. Finally, if $a\in \mathfrak r(M)=\mathfrak r(I)$, then $a^N\in I$
for some $N\geq 1$. Hence, if $b\in \cA$, then $ba^m\in I\subseteq M$ for all $m\geq N$ (since $I$ is an ideal). So $M$ is an MZ-space of $\cA$.
\epfv

\begin{lemma}\label{lma1.5} Let $\cA$ be a commutative ring contained in a $\Q$-algebra and $a_1,\dots,a_n\in \cA$.
If there exists $r\geq 0$ such that 
$$a_1^{r+i}+a_2^{r+i}+\cdots+a_n^{r+i}=0$$
\noindent for all $1\leq i\leq n$, then
all $a_j$ are nilpotent.
\end{lemma}

\noindent{\bf Proof:} If $\cA$ is an integral domain this is well-known and follows using a Vandermonde matrix. To show the general case we may assume that
$\cA$ is a $\Q$-algebra.
Then for each prime ideal $\mathfrak{p}$ of $\cA$, the quotient ring $\cA/\mathfrak{p}$ is an integral domain containing $\Q$. We obtain
from the integral domain case that each $a_j$ belongs to $\mathfrak{p}$ and hence to the intersection of all prime ideals of $\cA$, which is 
equal to the nil-radical of $\cA$. So each $a_j$ is nilpotent.
\epfv

\section{\bf The LFED Conjecture on Local Integral Domains}

Throughout this section $\cA$ is an {\bf integral domain} containing a field $k$ of characteristic zero, $D$ is a $k$-derivation on $\cA$ and $\f$ is a $k$-algebra endomorphism of $\cA$. By $\cA^*$ we denote the set of units of $\cA$ and by $Q(\cA)$ the quotient field of $\cA$.
We start with the following two useful observations:  

\begin{enumerate}
  \item[$i)$] \label{Obser1} {\it If $u\in \cA^*$ and $D$ is locally nilpotent, then $D(u)=0$}: note first $D(u)=au$ with $a\!:=u^{-1}D(u)$, and then apply Proposition \ref{Prop1.1}.  
  \item[$ii)$] {\it If $\cA$ is a local ring with maximal ideal $\mathfrak{m}$ and $a\in \cA$, then 
  $1-a\notin\mathfrak{m}$ or $1-2a\notin\mathfrak{m}$}: namely if both belong to $\mathfrak{m}$, then $a=(1-a)-(1-2a)\in\mathfrak{m}$, hence $1-a\notin\mathfrak{m}$, contradiction.
\end{enumerate}

\begin{lemma}\label{lma2.1} 
If $u\in \cA$ is transcendental over $k$ and $a_1,a_2,\dots\in k$ are distinct, then
the elements $\frac{1}{u-a_i}$ $(i\geq 1)$ in $Q(\cA)$ are linearly independent 
over $k$.
\end{lemma}

\noindent{\bf Proof:} Let $d\geq 1$ and $c_1,\dots,c_d\in k$ be such that $\sum_{i=1}^d c_i\frac{1}{u-a_i}=0$. By clearing up the denominators we see that $u$ is a root in $\cA$ of the polynomial 
$f(t)\!:=\sum_{i=1}^d c_ih_i(t)\in k[t]$, where $h_i(t)$ is the product of $t-a_j$ $(1\le j\ne i\le d)$. Since $u$ is transcendental over $k$, we have $f(t)=0$. 
Hence $f(a_i)=0$ all $1\le i\le d$. Then by the fact that $\cA$ is an integral domain 
we get $c_i=0$ for all $1\le i\le d$.   
\epfv

\begin{propo}\label{Prop2.2} 
If $\cA$ is a local ring and $D$ is locally finite, then $D=0$. 
\end{propo}

\noindent{\bf Proof:}  Let $a\in \cA$ and $V\!:=\sum_{j\geq 0} kD^j(a)$. We will show  that all eigenvalues of $D|_{V}$ are zero.
This implies that $D$ is locally nilpotent on $\cA$. By the observation $ii)$ above $1-a\notin\mathfrak{m}$ or $1-2a\notin\mathfrak{m}$. So 
$1-a$ or $1-2a$ belongs to $\cA^*$. By the observation $i)$ above we get that $D(a)=0$ or $D(2a)=0$. Since $\Q\subseteq \cA$ we get $D(a)=0$.
So $D=0$, as desired. 

Now, let $ v\in V$ with $D(v)=cv$ for some $c\in k$. Note also that $D(2v)=c(2v)$. So by the observation $ii)$ above we may assume that $1-v\notin\mathfrak{m}$, hence $(1-v)\in \cA^*$. If $v$ is algebraic over $k$, then $D(v)=0$ (for if $p(v)=0$, with $p(t)\in k[t]$ of the least degree, then $p'(v)D(v)=0$ implies $D(v)=0$) and we are done. So we may assume
that $v$ is transcendental over $k$. By induction on $m$ one easily shows that for all $m\geq 1$ we have 
\begin{align}\label{Prop2.2-peq1}
D^m((1-v)^{-1})=m!c^mv^m(1-v)^{-(m+1)}+(1-v)^{-m}p_m(v)
\end{align}
\noindent for some polynomial $p_m(v)\in k[v]$. Since $D$ is locally finite, there exist $N\geq 1$ and $c_1,\dots,c_N\in k$, with $c_N=1$, 
such that
$$
\sum_{m=1}^N c_mD^m((1-v)^{-1})=0.
$$
\noindent Substituting the formulas in Eq.\,(\ref{Prop2.2-peq1}) for the $D^m((1-v)^{-1})$ in the equation above we find, after multiplying 
by $(1-v)^{N+1}$, that $1-v$ divides $N!c^Nv^N$ (everything here is inside the polynomial ring $k[v]$). Then by the assumption that $v$ is transcendental over $k$ it is easy to see $c=0$. \epfv

\begin{propo}\label{Prop2.4} Assume that $\cA$ is a local ring and $\f$ is locally finite. If $\f(v)=cv$, with $0\neq v\in \cA$ and $0\neq c\in k$, then
\begin{enumerate}
  \item[$i)$] $c$ is a root of unity.
  \item[$ii)$] There does not exist $w\in \cA$, linearly independent from $v$, such that $\f(w)=cw+v$. 
\end{enumerate}
\end{propo}

\noindent{\bf Proof:} $i)$ If $v$ is algebraic over $k$, there exists $0\neq p(t)\in k[t]$ of the least degree such that $p(v)=0$. Applying $\f^m$ $(m\geq 1)$ 
 gives that $p(c^m v)=0$. So there exists $d\geq 1$ such that $c^d=1$. Now assume $v$ is transcendental over $k$. As in the proof of Proposition \ref{Prop2.2} we may assume that $(1-v)^{-1}\in \cA$. Since $\f$ is locally finite the elements
$\{\f^m(1/1-v)\}_{m\geq 1}$ are linearly dependent over $k$. Put $b\!:=1/v\in Q(\cA)$. So the elements 
$$
\frac{1}{1-c^mv}=b\,\frac{1}{b-c^m}
$$
\noindent are linearly dependent over $k$, and hence so are the elements 
$\{\frac{1}{b-c^m}\}_{m\geq 1}$.
If all $c^m$ are distinct, this contradicts lemma \ref{lma2.1} (applied to $Q(\cA))$. 
So $c^i=c^j$ for some $i<j$, which implies $i)$, since $c\neq 0$. 

$ii)$ Suppose there exists $w\in\cA$ such that $\f(w)=cw+v$. As above we may assume that 
$(1-w)^{-1} \in \cA$. Note that  $\f^m(w)=c^mw+mc^{m-1}v$ for all $m\geq 1$, 
and by $i)$ $c^d=1$ for some $d\ge 1$. 
So for all $n\geq 1$ we get
$$
\f^{nd+1}(\frac{1}{1-w})=\frac{1}{1-(cw+(nd+1)v)}=\frac{1}{v}\frac{1}{u-(nd+1)} 
$$
\noindent with $u=(1-cw)/v$.  Since $\f$ is locally finite
 these elements are linearly dependent over $k$. To get a contradiction we shall show that these elements are linearly independent over $k$. To do so, by lemma \ref{lma2.1} it suffices to show that $u$ is transcendental over $k$. 

Now, assume that $p(u)=0$ for some $0\neq p(t)\in k[t]$.  
Applying $\f^{nd}$ $(n\ge 1)$ to the equation $p(u)=0$ gives that
 $p(\f^{nd}(u))=0$ for all $n\geq 1$. So we have $\f^{n_1d}(u)=\f^{n_2d}(u)$ for some $1\leq n_1<n_2$. Note that for all $n\ge 1$ we have
 $$
\f^{nd}(u)=(1-c(w+nd c^{-1}v))/v=(1-cw-ndv)/v.
$$ 
The equality $\f^{n_1d}(u)=\f^{n_2d}(u)$ gives that $v=0$, a contradiction.
 \epfv
 
%
%

As a consequence of Proposition \ref{Prop2.4} we get the following: 

\begin{theo}\label{Thm2.4} 
Assume that $k$ is algebraically closed,
$\cA$ is a local ring and $\f$ is locally finite over $k$.
If $\cA$ is finitely generated as a $k$-algebra, or $\f$ is injective and $Q(\cA)$ 
is a finitely generated field extension of $k$, then
there exist $1\leq i<j$ such that $\f^i=\f^j$. In particular, 
$\im (1-\f)$ is an MZ-space of $\cA$.
\end{theo}

\noindent{\bf Proof:} 
First, assume that $\cA=k[a_1,\dots,a_n]$ for some $a_i\in \cA$. 
Put $W\!:=\sum_{i=1}^n k a_i$ and $V\!:=\sum_{j\geq 0}\f^j(W)$. Since $\f$ is locally finite, $V$
is finite dimensional. By Proposition \ref{Prop2.4} each nonzero eigenvalue 
of $\f|_V$ is a root of unity and each Jordan block of the representative matrix of $\f|_V$ 
(with respect to any fixed $k$-linear basis of $V$) associated with a nonzero eigenvalue 
is an $1\times 1$ block. Then by the Jordan-Chevalley decomposition of $\f$ 
(e.g., see \cite[Proposition $1.3.8$]{E}) and   
 \cite[Proposition $4.2$]{H}) we have $V=V_0\oplus V_+$, where
$\f^N(V_0)=0$ for some $N\geq 1$ and $V_+=\oplus_{j=1}^r k v_j$ for some $r\geq 0$ and $v_j\in \cA$ with 
$\f(v_j)=c_j v_j$ and $c_j$ being 
a root of unity. Let $d\geq 1$ be such that $c_j^d=1$ for all $j$. 
Finally choose a $k$-basis $u_1,\dots, u_s$ of $V_0$.
Then obviously $\cA=k[v_1,\dots, v_r,u_1,\dots,u_s]$. Since $\f^d(v_j)=v_j$ for all $j$ and $\f^N(u_i)=0$ for all $i$ it follows that $\f^{N+d}=\f^N$ on $\cA$. Consequently, $\im (1-\f)$ is an MZ-space of $\cA$ by 
Proposition \ref{Prop1.4}.

Now, assume that $Q(\cA)=k(a_1, \dots, a_n)$ for some $a_i\in \cA$. 
Define $W$, $V$, $V_0$ and $v_i$'s as above. Observe
that $V_0=0$ since $\f$ is injective. Then $Q(\cA)=k(v_1,\dots,v_r)$. Since $\f$ is injective we can extend $\f$ to $Q(\cA)$. Since $\f^d(v_j)=v_j$ for all $j$ we deduce that $\f^d=1$ on $Q(\cA)$ 
and hence also on $\cA$. Then by Proposition \ref{Prop1.4} the theorem follows again.
\epfv


Let $x=(x_1, x_2, \dots, x_n)$ be $n$ free commutative variables and $k(x)$ the quotient field of the polynomial algebra $k[x]$. Then $k(x)$ is a field (hence, a local integral domain), and each $k$-algebra endomorphism $\f$ of $k(x)$ is injective. By Proposition \ref{Prop2.2} and Theorem \ref{Thm2.4} we have the following:

\begin{corol}\label{k(x)Thm}
The LFED Conjecture holds for $k(x)$.
\end{corol} 

\noindent{\bf Proof:} 
Let $\bar k$ be the algebraic closure of $k$. Then it is well-known that $\bar k\otimes_k k(x)\simeq\bar k(x)$. So by the remark
made after the introduction of the LFED Conjecture we may assume that $k$ is algebraically closed.
Then the corollary follows immediately from Proposition \ref{Prop2.2} and Theorem \ref{Thm2.4}. 
\epfv

Next we show the LFED Conjecture for the formal power series algebra $k[[x]]$ (in $n$ variables), which is a local integral domain but is not finitely generated as a $k$-algebra.

\begin{propo}\label{Prop2.7} Let $\f:k[[x]]\rightarrow k[[x]]$ be a locally finite $k$-algebra endomorphism of $k[[x]]$. Then there exist $1\leq i<j$ such that $\f^i=\f^j$. In particular, $\im(1-\f)$ is a MZ-space of $k[[x]]$.
\end{propo}

\noindent{\bf Proof:}  Let $\bar{k}$ be the algebraic closure of $k$.  Replacing $k[[x]]$ by $\bar{k}[[x]]=\bar{k}\otimes_k k[[x]]$, 
again using the remark made after the introduction of the LFED Conjecture, we may assume that $k$ is algebraically closed. 
Note also that the last statement of the proposition follows from the first one and Proposition \ref{Prop1.4}.

Put $W\!:=\sum_{i=1}^n k x_i$ and $V\!:=\sum_{j\geq 0}\f^j(W)$. Since $\f$ is locally finite, $V$ is finite dimensional. Applying the Jordan-Chevalley decomposition of $\f$ 
and Proposition \ref{Prop2.4} (similarly as in the proof of Theorem \ref{Thm2.4}) we have $V=V_0\oplus V_+$, where
$\f^N(V_0)=0$ for some $N\geq 1$ and $V_+=\oplus_{j=1}^r k v_j$ for some $r\geq 0$ with $\f(v_j)=c_j v_j$ and $c_j$
is a root of unity. Let $d\geq 1$ be such that $c_j^d=1$ for all $j$. Then it follows that $\f^{N+d}=\f^d$ on $V$, hence
on each $x_i$. 

To show that $\f^{N+d}=\f^d$ on $k[[x]]$ it suffices to prove that $\f$ is continuous with respect to
the $\mathfrak{m}\!:=(x_1,\dots,x_n)$-adic topology on $k[[x]]$. It suffices to show that $\f(x_i)\in\mathfrak{m}$ for all $i$.
Suppose that $\f(x_i)=c+h_i(x)$, with $c\in k^*$ and $h_i(x)\in\mathfrak{m}$.  Then $(x_i-c)^{-1}\in k[[x]]$, hence $\f((x_i-c)^{-1})\in k[[x]]$, i.e., $1/((c+h_i(x))-c)=1/h_i(x) \in k[[x]]$, a contradiction 
since $h_i(0)=0$.
\epfv

The above argument can be sharpened to obtain a description for all
locally finite $k$-endomorphisms of $k[[x]]$, in case $k=\bar{k}$: 

\begin{theo}\label{Thm2.8} Assume that $k=\bar k$ and let $\f$ be 
a locally finite endomorphism of $k[[x]]$. Then there exist $N\ge 1$, 
$1\leq d\leq n$, $c_1,\dots,c_d$ roots of unity in $k$ and coordinates
$y_1, \dots, y_n\in k[[x]]$, i.e., $k[[x]]=k[[y_1,\dots,y_n]]$, such that $\f(y_i)=c_iy_i$ for all $1 \le i\leq d$ and $\f^N(y_i)=0$ for all $i>d$ (if $d=n$ the last statement is empty).
\end{theo}

\noindent{\bf Proof:} We use the same notations fixed in the proof of Proposition \ref{Prop2.7}.
 For $v\in V$ let $v_1$ denote the linear part of $v$. Choose a $k$-basis $u_1,\dots, u_t$ of $V_0$. It follows that 
$$
W\subseteq\sum k v_{j1}+\sum k u_{i1}\subseteq\sum k x_i=W.
$$
\noindent So the elements $v_{j1}$ and $u_{i1}$ form a spanning set of $W$. Hence there exist $u_{i_11},\dots,u_{i_d 1}$, $v_{j_11}, \dots, v_{j_{n-d}1}$ which form a $k$-basis of $W$ and whose determinant of the Jacobian matrix satisfies 
$$
\det\,J(u_{i_1},\dots,u_{i_d},v_{j_1},\dots,u_{j_{n-d}})(0)\in k^*.
$$
\noindent So by the formal inverse function theorem we obtain that
$$
k[[x]]=k[[u_{i_1},\dots,u_{i_d}, v_{j_1},\dots,u_{j_{n-d}}]],
$$ 
\noindent from which and Proposition \ref{Prop2.4} the theorem follows.
\epfv

By Propositions \ref{Prop2.2} and \ref{Prop2.7} we immediately have the following:

\begin{corol}\label{k[[x]]Thm}
The LFED Conjecture holds for $k[[x]]$.
\end{corol}

\section{\bf The LFED Conjecture on $\mathbf{k[[x]][x^{-1}]}$}

Recall that the Duistermaat-van der Kallen theorem follows if some special cases of the LFED Conjecture holds for the algebra $\C[x,x^{-1}]$. On the other hand the Duistermaat-van der Kallen theorem can not be generalized to the algebra $\C[[x]][x^{-1}]$ (in the obvious way). To see this take $a=1/x_1$
and $b=1/(1-x_1)$. Then clearly all powers of $a$ have no constant term, however all elements $ba^m$  
($m\geq 1$) have a non-zero constant term. It is therefore natural to ask if the LFED Conjecture also fails for the algebra $\C[[x]][x^{-1}]$.

The aim of this section is to show that it does not. Throughout this section we have the following notations: $k$ is a field of characteristic zero; $k[[x]]\!:=k[[x_1,\dots,x_n]]$ is the algebra of the formal power series over $k$; and finally $k[[x]][x^{-1}]$ denotes the localization of $k[[x]]$ in the elements $x_i$ $(1\le i\le n)$. 

We first show that the derivation case of the LFED Conjecture holds for $k[[x]][x^{-1}]$. This follows immediately from the following: 

\begin{theo}\label{Thm2.5} 
$k[[x]][x^{-1}]$ has no nonzero locally finite $k$-derivations.
\end{theo}
 
\noindent{\bf Proof:} Let $D$ be a locally finite $k$-derivation on $\cA\!:=k[[x]][x^{-1}]$. Assume $D\neq 0$. Write $D=\sum_{i=1}^n p_i\p_i$,
with $p_i\in \cA$ for all $i$. We will show that all $p_i\in k[[x]]$. From this it follows that the restriction of $D$ to $k[[x]]$ is
a locally finite derivation on the local algebra $k[[x]]$. Then by Proposition \ref{Prop2.2} we have 
$D=0$ on $k[[x]]$ and hence also on $\cA$.

Suppose that $p_i\notin k[[x]]$ for some $i$. Then there exists a monomial $x^a$ in $p_i$, with a non-zero coefficient, such 
that $a_j<0$ for some $j$, where $a=(a_1,\dots,a_n)\in\Z^n$. Choose $d_j\in\N$ large enough such that
$d_ja_j+\sum_{i\neq j} a_i<0$ and put $d\!:=(1,\dots,1,d_j,1,\dots,1)$. Then $\cA$ becomes a $\Z$-graded algebra by defining the weight of a monomial $x^b$ by $w(x^b)\!:=\langle d,b\rangle$, where $\langle,\rangle$ denotes the usual inner-product on $\R^n$.
So $w(x^a)<0$ and $\cA=\oplus_{r\in\Z} \cA_r$, where $\cA_r$ is the $k$-span of all monomials having weight $r$. Furthermore, we
can write $D$ uniquely as $D=\sum_{i\geq N} D_i$, with $N\in\Z$, $D_N\neq 0$ and $D_i$ a $k$-derivation on $k[x,x^{-1}]$ 
such that $D_i\cA_r\subseteq \cA_{r+i}$ for all $r\in\Z$. Since $w(x^a)<0$ it follows that $N<0$. 

Now fix an arbitrary $b\in\Z^n$. The fact that $N<0$ and $D$ is locally finite imply that $D_N^m(x^b)=0$ for all large $m$. Consequently, $D_N$ is locally nilpotent on $k[x,x^{-1}]$. Since each $x_s$ is a unit in $k[x,x^{-1}]$, it follows from observation $i)$,  Section $2$ that $D(x_s)=0$ for all $1\leq s\leq n$. 
So $D=0$, a contradiction.
\epfv

Next we  show the $\mathcal{E}$-derivation case of the LFED Conjecture for the algebra 
$k[[x]][x^{-1}]$ in the following:

\begin{theo}\label{Thm2.6}  Let  $\cA\!:=k[[x]][x^{-1}]$. If $\f:\cA\rightarrow \cA$ is a locally finite $k$-algebra homomorphism, then there exist
$1\leq i<j$ such that $\f^i=\f^j$. In particular, $\im (1-\f)$ is an MZ-space of $\cA$. 
\end{theo}

To prove this theorem we need to show first the following: 
 
\begin{lemma}\label{Lma2.9} 
Let $\cA\!:=k[[x]][x^{-1}]$ and $\f:\cA\rightarrow \cA$ be a $k$-algebra homomorphism (not necessarily locally finite). Then $\f(k[[x]])\subseteq k[[x]]$.
\end{lemma}

\noindent{\bf Proof:} As before, we may again assume that $k$ is algebraically closed. We will show that $\f(x_i)\in k[[x]]$ for all $i$. Suppose not, say, $\f(x_1)\notin k[[x]]$. Since $\f(x_1)$
is a unit in $\cA$ (for $x_1$ is in $\cA$), there exist $\al\in\Z^n$ and $h\in k[[x]]$ with $h(0)\neq 0$ such that $\f(x_1)=x^{\al}h$.
Since $\al\notin\N^n$ we can write $\al=a-b$, with $a,b\in\N^n$ such that $b\neq 0$ and for each $i$ either $a_i=0$
or $b_i=0$. So we may write $x^{\al}=x^a/x^b$. Then 
$$
\f(1+x_1)=1+\f(x_1)=1+\frac{x^ah}{x^b}=x^{-b}(x^b+x^ah).
$$

Since $\f(1+x_1)$ is a unit in $\cA$ (for $1+x_1$ is a unit in $k[[x]]$), it follows from the choices of $a$ and $b$ that $a=0$ (since $b\neq 0$). So $\f(1+x_1)= x^{-b}(x^b+h)$. Now let $r\geq 1$. Observe that $1+x_1=u_r(x)^r$
for some $u_r(x)\in k[[x]]$ with $u_r(0)=1$. So $\f(u_r)$ is a unit in $\cA$ and hence also 
of the form $x^{c_r}g_r(x)$, with 
$c_r\in\Z^n$ and $g_r(x)\in k[[x]]$ with $g_r(0)\neq 0$. Consequently, 
$$
x^{-b}(x^b+h)=\f(1+x_1)=x^{rc_r}g_r(x)^r.
$$
\noindent Since $h(0)\neq 0$ and $g_r(0)\neq 0$, we get that $-b=rc_r$ for all $r\geq 1$. 
Hence $b=0$, a contradiction. 
\epfv

Finally, we are able to show Theorem \ref{Thm2.6} as follows. 

\medskip

\noindent{\bf Proof of Theorem \ref{Thm2.6}.} It follows from Lemma \ref{Lma2.9} that by restricting to $k[[x]]$ the map $\f$ induces
a $k$-algebra endomorphism of $k[[x]]$, which is obviously also locally finite.
Then Proposition \ref{Prop2.7} implies that there exist $1\leq i<j$ such that $\f^i=\f^j$ on $k[[x]]$. Hence the same holds on $\cA$. Then the theorem follows immediately from Proposition \ref{Prop1.4}.
\epfv

\medskip

\noindent {\bf Acknowledgement.} The authors would like to thank the referee for some useful suggestions and remarks.

\end{document}